\documentclass{article}

\usepackage{amsmath}
\usepackage[latin1]{inputenc}
\usepackage[T1]{fontenc}
\usepackage[english]{babel}
\usepackage{enumerate}
\usepackage{amssymb}
\usepackage{amsthm}
\usepackage{cleveref}
\usepackage{titlesec} 
\usepackage{url}

\titleformat{\section}[block] 
  {\normalfont\bfseries\centering} 
  {} 
  {0pt} 
  {}

\newtheorem{Lemma}{Lemma}[section]
\newtheorem{Theorem}[Lemma]{Theorem}

\title{A Finite Dimensional Counterexample for Arveson's Hyperrigidity Conjecture}

\date{}
\author{Marcel Scherer}

\begin{document}

\maketitle
\begin{abstract}
We construct an operator system generated by $4$ operators that is not hyperrigid, although all restrictions of irreducible representations have the unique extension property. 
\end{abstract}

\section{Introduction}
\let\thefootnote\relax\footnotetext{ \date\\
2020 \textit{Mathematics Subject Classification.}Primary 46L07, 47A20; Secondary: 46L52.\\
\ \textit{Key words and phrases.} operator system, hyperrigidity conjecture, unital completely positive, affine function.\\
 The author was partially supported by the Emmy Noether Program of the German Research Foundation (DFG Grant 466012782).}

Let $S$ be a separable operator system and $C^*(S)$ a $C^*$-algebra generated by $S$. Arveson conjectured in \cite{WA} that every representation of $C^*(S)$, restricted to $S$, has a unique extension to a unital completely positive map if and only if every irreducible representation of $C^*(S)$, restricted to $S$, has a unique unital completely positive extension - a conjecture known as Arveson's hyperrigidity conjecture. For a detailed background, we recommend \cite{BD}.\\
It has been shown by Bilich and Dor-On in \cite{BD} that the conjecture fails for an infinite-dimensional operator system generated by an operator algebra in a non-commutative $C^*$-algebra. However, a critical part of that proof relies on the infinite dimensionality of the operator system. In this paper, we construct a finite-dimensional counterexample for Arveson's hyperrigidity conjecture in \Cref{Theorem}.\\
We use the main ideas of \cite{BD}, with the only significant difference being in the proof of the unique extension property for the irreducible representations of $C^*(S)$. For this, we employ a technique found in \cite{AW}, which can be used to show that a certain projection is below a family of positive multiplication operators, and that pure states have a maximal irreducible dilation \cite{KK}.\\

\section{A Finite Dimensional Counterexample}
Let $\mathbb{B}_4$ be the Euclidean open unit ball in $\mathbb{R}^4$, and let $S_3$ be its boundary, the sphere in $\mathbb{R}^4$. Let $A(\overline{\mathbb{B}_4})$ denote the continuous affine functions on the closed unit ball, and let $t_i\in A(\overline{\mathbb{B}_4})$, $i=1,2,3,4$, be the projection onto the $i$-th component. Note that $A(\overline{\mathbb{B}_4})$ is spanned by $1, t_1, t_2, t_3, t_4$.\\
Define $P:  L^2(S_3)\to \mathbb{C}, g\mapsto \langle g, 1\rangle$. For $i=1,2,3,4$ and $c>0$, we define the operators 
  \begin{equation*}
      T_i=M_{t_i}\oplus 0,\ T_{1,c}=\begin{pmatrix} M_{t_1} & cP^* \\ cP & 0\end{pmatrix},\ \tilde T_{1,c}=\begin{pmatrix} M_{t_1} & cP^* \\ 0 & 0\end{pmatrix}
  \end{equation*}
on $L^2(S_3)\oplus \mathbb{C}$. Here, $L^2(S_3)$ is equipped with the unique rotation-invariant probability measure $m$ on $S_3$. Let $S_c$ be the operator system generated by $\tilde T_{1,c}, T_2, T_3, T_4$.  In \Cref{Theorem}, we show that $S_c$ is not hyperrigid, but the restrictions of all irreducible representations of $C^*(S_c)$ are boundary representations.\\ 
A crucial part of the proof is to show that the \textit{joint numerical range} of the operator tuple $(T_{1,c},T_2,T_3,T_4)$, defined as
  \begin{equation*}
    \mathcal{W}((T_{1,c},T_2,T_3,T_4))=\{(\langle T_{1,c}x,x\rangle,\langle T_2x,x\rangle,\langle T_3x,x\rangle,\langle T_4x,x\rangle);\ \substack{x\in L^2(S_3)\oplus\mathbb{C},\\ \|x\|=1}\},
  \end{equation*} 
is contained in $\mathbb{B}_4$. To prove this, we first need to show that $\mathcal{W}((T_1, T_2, T_3, T_4))$ is contained in $\mathbb{B}_4$.

\begin{Lemma}
It holds that:
  \begin{equation*}
    \mathcal{W}((T_1,T_2,T_3,T_4))\subset\mathbb{B}_4.
  \end{equation*}
\label{Lemma0}
\end{Lemma}

\underline{Proof}:\\
Let $x\in  L^2(S_3)\oplus\mathbb{C}$ with $\|x\|=1$. Using the fact that $\sum_{i=1}^4t_i^2=1$ on $S_3$ and applying the Cauchy-Schwarz inequality, we get:
  \begin{equation*}
    \|(\langle T_1x,x\rangle,\langle T_2x,x\rangle,\langle T_3x,x\rangle,\langle T_4x,x\rangle)\|^2=\sum_{i=1}^4(\langle T_ix,x\rangle)^2\le\sum_{i=1}^4\langle T_i^2x,x\rangle=\|x\|
=1.
  \end{equation*}
Furthermore, the Cauchy-Schwarz inequality implies that if equality holds, there must exists some $i\in\{1,2,3,4\}$ and a $0\neq\lambda\in\mathbb{C}$ such that $T_ix=\lambda x$. However, $M_{t_i}$ has no eigenvalues, so we conclude that:  
  \begin{equation*}
    \sum_{i=1}^4(\langle T_ix,x\rangle)^2<1,
  \end{equation*} 
which completes the proof. $\hfill\square$\\

\begin{Lemma}
Let $0<c<1/2$. Then, the following holds:
  \begin{equation*}
    \mathcal{W}((T_{1,c},T_2,T_3,T_4))\subset \mathbb{B}_4.
  \end{equation*}
\label{Lemma1}
\end{Lemma}
\underline{Proof}:\\
Let $0<c\le1/2,$ and let $\tilde S_c$ be the operator system generated by $T_{1,c}, T_2, T_3, T_4$. Define the map
  \begin{equation*}
    \Phi: A(\overline{\mathbb{B}_4})\to \tilde S_c, \alpha+\beta t_1+\gamma t_2+\delta t_3+\epsilon t_4\mapsto \alpha+\beta T_{1,c}+\gamma T_2+\delta T_3+\epsilon T_4.
  \end{equation*}
It is clear that this map is well-defined, bijective and that its inverse is positive. The next step is to show that $\Phi$ is positive. So let $0\le\alpha+\beta t_1+\gamma t_2+\delta t_e+\epsilon t_4=f$ and notice that this implies $\alpha\ge0$ and $\beta, \gamma, \delta, \epsilon\in\mathbb{R}$. If $\alpha=0$, then $f=0$, and there is nothing to show. Thus, assume $\alpha>0$. Then:
  \begin{equation*}
    \Phi(f)=\begin{pmatrix} M_f & c\beta P^* \\ c\beta P & \alpha \end{pmatrix}.
  \end{equation*}    
is positive if and only if the Schur complement 
  \begin{equation*}
    M_f-c^2\beta^2\alpha^{-1}P^*P
  \end{equation*}
is positive (see, for example, \cite[Lemma 7.2.7]{AM}). Therefore, the positivity of $\Phi(f)$ is equivalent to
  \begin{equation*}
    c^2\beta^2 P^*P\le \alpha M_f.
  \end{equation*}
Evaluating $f$ in $(1,0,0,0)$ and $(-1,0,0,0)$ shows that $|\beta|\le\alpha$, and since there is nothing to show for $\beta=0$, it suffices to check that
  \begin{equation*}
    c^2P^*P\le|\beta|^{-1}M_f.
  \end{equation*}
Recall that $m$ is the unique rotation-invariant probability measure on $S_3$. Let $g\in L^2(S_3)$, and write $\tilde f=f/|\beta|$. Following the idea of \cite{AW}, we get
  \begin{equation*}
    \begin{split}
      \langle P^*Pg,g\rangle&=\left|\int_{S_3}g dm\right|^2\\
      &\le\left(\int_{S_3} \tilde f^{1/2}|g|\tilde f^{-1/2}dm\right)^2\\
      &\le\left(\int_{S_3}\tilde f|g|^2dm\right)\left(\int_{S_3}\tilde  f^{-1}dm\right)\\
      &=\langle M_{\tilde f}g,g\rangle\left(\int_{S_3} \tilde f^{-1}dm\right).
    \end{split}
  \end{equation*}
Note that $\tilde f(z)=|\beta|^{-1}(\alpha + \langle z, \omega\rangle_{\mathbb{R}^4})$, where $\omega=(\beta,\gamma,\delta,\epsilon)$, and $\alpha\ge\|\omega\|$ since $f\ge0$. Let $U\in\mathcal{B}(\mathbb{R}^4)$ be an orthogonal matrix such that $U^*\omega=(\|\omega\|,0,0,0)$. Then
  \begin{equation*}
    \begin{split}
      \int_{S_3}\tilde f^{-1}dm&=\int_{S_3}\frac{|\beta|}{\alpha+\langle z,\omega\rangle}dm(z)=\int_{S_3}\frac{|\beta|}{\alpha+\langle U(z),\omega\rangle}dm(z)\\
      &=\int _{S_3}\frac{|\beta|}{\alpha+\|\omega\|t_1}dm\le\int_{S_3}\frac{|\beta|}{\|\omega\|}\frac{1}{1+t_1}dm\le\int_{S_3}\frac{1}{1+t_1}dm
    \end{split}
  \end{equation*}
and 
  \begin{equation*}
\int_{S_3}\frac{1}{1+t_1}dm=(2\pi^2)^{-1}\int_{0}^{\pi}\int_{0}^{\pi}\int_{0}^{2\pi}\frac{\sin^2(x)\sin(y)}{1+\cos(x)}dzdydx=2.
  \end{equation*}
Therefore, $\Phi(f)\ge0$ if $c\le 1/2$. Thus, we have shown that if $\phi$ is a positive state on $\tilde S_c$, then $\phi\circ\Phi$ is a positive state on $A(\overline{\mathbb{B}_4})$, and therefore 
  \begin{equation*}
    (\phi(T_1), \phi(T_2), \phi(T_3), \phi(T_4))\in\overline{\mathbb{B}_4}.
  \end{equation*}
 In particular, for $0<c\le1/2$, we obtain
  \begin{equation}
    \mathcal{W}((T_{1,c},T_2,T_3,T_4))\subset\overline{\mathbb{B}_4}.
\label{Gl1}
  \end{equation}
Let $x\in L^2(S_3)\oplus\mathbb{C}$ with $\|x\|=1$ and define
  \begin{equation*}
    \begin{split}
      &z_1=(\langle T_1 x,x\rangle, \langle T_2x,x\rangle, \langle T_3x,x\rangle, \langle T_4x,x\rangle),\\
      &z_2=(\langle T_{1,1/2} x,x\rangle, \langle T_2x,x\rangle, \langle T_3x,x\rangle, \langle T_4x,x\rangle).
    \end{split}
  \end{equation*}
Then, $z_1\in\mathbb{B}_4$ by \Cref{Lemma0} and $z_2\in\overline{\mathbb{B}_4}$ by \Cref{Gl1}. Thus, for $0<c<1/2$, we have that 
  \begin{equation*}
    (\langle T_{1,c} x,x\rangle, \langle T_2x,x\rangle, \langle T_3x,x\rangle, \langle T_4x,x\rangle)=(1-2c)z_1+2cz_2\in\mathbb{B}_{4}.
\label{Eq1}
  \end{equation*}
$\hfill\square$\\
\ \\
\begin{Lemma}
Let $0<c<1$ and $x\in L^2(S_3)\oplus\mathbb{C}$ with $\|x\|=1$. Then, the following holds:
  \begin{equation*}
    \|(\langle \tilde T_{1,c}x,x\rangle, \langle T_2x,x\rangle, \langle T_3x,x\rangle, \langle T_4x,x\rangle)\|<1.
  \end{equation*}
\label{Lemma2}
\end{Lemma}

\underline{Proof}:\\
Let $x=(y,a)\in L^2(S_3)\oplus \mathbb{C}$ with $\|x\|=1$ and define
  \begin{equation*}
    z=(\langle \tilde T_{1,c}x,x\rangle, \langle T_2x,x\rangle, \langle T_3x,x\rangle, \langle T_4x,x\rangle).
  \end{equation*} 
Assume, for contradiction, that $\|z\|\ge1$. The equation
 \begin{equation*}
    |\langle \tilde T_{1,c}x,x\rangle|= |\langle M_{t_1}y,y\rangle+ca\langle 1,y\rangle|\le |\langle M_{t_1}y,y\rangle|+|ca\langle 1,y\rangle|
  \end{equation*}
shows, on one hand, that $\langle 1,y\rangle\neq 0$, and on the other, that for 
  \begin{equation*}
    \tilde x=\begin{cases} (y,|a|\frac{\langle y,1\rangle}{|\langle 1,y\rangle|}) &\textup{if} \ \langle M_{t_1}y,y\rangle\ge0\\(y,-|a|\frac{\langle y,1\rangle}{|\langle 1,y\rangle|}) & \textup{else} \end{cases}
  \end{equation*}
and
  \begin{equation*}
    \tilde z=(\langle \tilde T_{1,c}\tilde x,\tilde x\rangle, \langle T_2\tilde x,\tilde x\rangle, \langle T_3\tilde x,\tilde x\rangle, \langle T_4\tilde x,\tilde x\rangle),
  \end{equation*}
we have $\|\tilde x\|=\|x\|=1$, $\langle \tilde T_{1,c} \tilde x, \tilde x\rangle\in\mathbb{R}$ and $1\le\|z\|\le\|\tilde z\|$. Thus, applying the identity $\textup{Re}(\tilde T_{1,c})=T_{1,c/2}$, we conclude that 
  \begin{equation*}
    \tilde z=(\langle T_{1,c/2}\tilde x,\tilde x\rangle, \langle T_2\tilde x,\tilde x\rangle, \langle T_3\tilde x,\tilde x\rangle, \langle T_4\tilde x,\tilde x\rangle),
  \end{equation*}
which lies within the unit ball $\mathbb{B}_4$ by \Cref{Lemma1}, thus leading to a contradiction. Hence, $\|z\|<1$. $\hfill\square$\\

\ \\
\begin{Lemma}
Let $0<c<1$. Then, the compact operators $\mathcal{K}(L^2(S_3)\oplus\mathbb{C})$ are contained in $C^*(S_c)$.
\label{Lemma3}
\end{Lemma}

\underline{Proof}:\\
The proof is essentially the same as in \cite[Lemma 3.2]{BD}. We start by noting that
  \begin{equation*}
    \tilde T_{1,c}\tilde T^*_{1,c}+\sum_{i=2}^4T_i^*T_i-1=\begin{pmatrix} c^2P^*P & 0\\ 0 & -PP^* \end{pmatrix}\in C^*(S_c).
  \end{equation*}
This implies that
  \begin{equation*}
    0\oplus PP^*=\frac{1}{c^{-2}+1}\left(c^{-2}\begin{pmatrix} c^2P^*P & 0\\ 0 & -PP^* \end{pmatrix}^2-\begin{pmatrix} c^2P^*P & 0\\ 0 & -PP^* \end{pmatrix}\right)\in C^*(S_c) ,
  \end{equation*}
  \begin{equation*}
      id\oplus 0=1-0\oplus PP^*\in C^*(S_c).
  \end{equation*}
Therefore,
  \begin{equation*}
    P^*P\oplus 0=\frac{1}{c^2}(id\oplus 0)(\tilde T_{1,c}\tilde T^*_{1,c}+\sum_{i=2}^4T_i^*T_i-1)\in C^*(S_c).
  \end{equation*}
Additionally, $M_{t_1}\oplus0=\tilde T_{1,c}(id\oplus0)\in C^*(S_c)$ and
  \begin{equation*}
    \begin{pmatrix} 0 & 0\\ P & 0\end{pmatrix}=\tilde T_{1,c}^*-M_{t_1}\oplus0 \in C^*(S_c).
  \end{equation*}
 Further multiplying by $M_{t_i}\oplus0$, $i=1,2,3,4$, from the left and right to $P^*P\oplus 0$ and $\begin{pmatrix} 0 & 0\\ P & 0\end{pmatrix}$ shows that 
  \begin{equation*}
  \begin{pmatrix} p\langle\cdot, q \rangle& 0\\ 0 & 0\end{pmatrix},  \begin{pmatrix} 0 & 0\\ \langle\cdot,q\rangle & 0\end{pmatrix}\in C^*(S_c)
  \end{equation*}
for every $p, q \in \mathbb{C}[t_1,t_2,t_3,t_4]$. Finally, since the polynomials are dense in $L^2(S_3)$, any compact operator can be approximated by elements of $C^*(S_c)$, completing the proof. $\hfill\square$\\
\ \\
The previous lemma places us in a setting similar to the discussion following \cite[Lemma 3.2]{BD}. Since $(T_{1,c},T_2,T_3,T_4)$ is a compact perturbation of $(M_{t_1}\oplus 1, M_{t_2}\oplus 1, M_{t_3}\oplus 1, M_{t_4}\oplus 1)$, by the previous lemma, we have the following split short exact sequence:
  \begin{equation*}
    0\to\mathcal{K}(L^2(S_3\oplus\mathbb{C}))\to C^*(S_c)\to C(S_3)\to 0.
  \end{equation*}
Consequently, $C^*(S_c)$ is a type $I$ $C^*$-algebra by \cite[Theorem 1]{JG}, and the only irreducible representations of $C^*(S_c)$ are given by the identity representation and the evaluations $e_z$, defined by 
  \begin{equation*}
    C^*(S_c)\to \mathbb{C}, (T_{1,c}, T_2, T_3, T_4)\mapsto z
  \end{equation*}
for $z\in S_3$, see \cite[Theorem 1.3.4]{WA2} and \cite[Page 20, Corollary 2]{WA2}.\\
Additionally, in the proof of the following theorem, we use the facts that a state on $S_c$ is pure if and only if it is a extreme point of the state space $S(S_c)$ and that a unital completely positive map $\phi$ on $S_c$ is maximal if and only if it has the unique extension property.

\begin{Theorem}
Let $0<c<1$. The operator system $S_c$ is not hyperrigid. However, the restrictions of all irreducible representations of $C^*(S_c)$ to $S_c$ have the unique extension property.
\label{Theorem}
\end{Theorem}

\underline{Proof}:\\
To show that $S_c$ is not hyperrigid, we begin by considering the $*$-homomorphism
  \begin{equation*}
    C^*(S_c)\to\mathcal{B}(L^2(S_3)),(\tilde T_{1,c}, T_2, T_3, T_4)\mapsto (M_{t_1}, M_{t_2}, M_{t_3}, M_{t_4}).
  \end{equation*}
This homomorphism dilates non-trivially to the identity representation and thus the restriction to $S_c$ does not have the unique extension property, proving that $S_c$ is not hyperrigid. \\
It remains to show that the irreducible representations of $C^*(S_c)$ are boundary representations. By Arvson's boundary theorem, we observe that since
  \begin{equation*}
    1=\|\sum_{i=1}^4 T^*_iT_i\|<\|T^*_{1,c}T_{1,c}+\sum_{i=2}^4 T_i^*T_i\|,
  \end{equation*}
the identity representation of $C^*(S_c)$ is a boundary representation. Thus, it remains to verify that the point evaluations $e_z$, restricted to $S_c$, are maximal for all $z\in S_3$.\\
We begin by showing that   
  \begin{equation}
    \|(\phi(\tilde T_{1,c}), \phi(T_2), \phi(T_3), \phi(T_4))\|\le1
\label{Eq3}
  \end{equation}
for every $\phi\in S(S_c)$ and that 
\begin{equation}
    \|(\phi(\tilde T_{1,c}), \phi(T_2), \phi(T_3), \phi(T_4))\|<1
\label{Eq4}
  \end{equation}
for every pure state $\phi$ that is not maximal.\\
Let $\phi\in S(S_c)$. If $\phi$ is pure and maximal, we have
  \begin{equation*}
    \|(\phi(\tilde T_{1,c}), \phi(T_2), \phi(T_3), \phi(T_4))\|=1,
  \end{equation*}
 since the only representations of $C^*(S_c)$ with image in $\mathbb{C}$ are given by the point evaluations. If $\phi$ is pure and not maximal, then by \cite[Theorem 2.4]{KK}, $\phi$ dilates non-trivially to a maximal irreducible u.c.p. map, which must be the identity representation, since this is the only irreducible representation besides the point evaluations. Thus, there exists $x\in L^2(S_3)\oplus \mathbb{C}$ with $\|x\|=1$ such that $\phi(\cdot)=\langle \cdot x,x\rangle$. Since $c<1$, we can apply  \Cref{Lemma2} to obtain \Cref{Eq4}, and since the convex hull of the extreme points of $S(S_c)$ is $S(S_c)$, by Carath\'edory's theorem, we also obtain \Cref{Eq3}.\\
It follows from \Cref{Eq3} that the restrictions of the maps $e_z$ to $S_c$ are extreme points of $S(S_c)$. By equation \Cref{Eq4}, these maps must also be maximal, completing the proof.$\hfill\square$

\section*{Acknowledgement}
The author thanks Michael Hartz for pointing out that a counter example could be constructed using the method in \cite{AW}.

\bibliography{bib_hyper_counterexample} 
\bibliographystyle{plain}

Fachrichtung Mathematik, Universit\"at des Saarlandes, 66123 Saarbr\"ucken, Germany\\
\textit{Email address:} scherer@math.uni-sb.de

\end{document}